\documentclass{jag-l}
\usepackage[arrow,matrix,graph]{xy}

\newtheorem{theorem}{Theorem}[section]
\newtheorem{corollary}[theorem]{Corollary}
\newtheorem{lemma}[theorem]{Lemma}

\theoremstyle{definition}

\newtheorem{remark}[theorem]{Remark}
\newtheorem*{ack}{Acknowledgments}
\newtheorem*{claim}{Claim}

\newcommand{\Z}{\mathbb{Z}}

\newcommand{\C}{\mathbb{C}}

\newcommand{\RR}{{\mathcal R}}

\newcommand{\G}{\Gamma}
\newcommand{\E}{\mathsf{E}}
\newcommand{\V}{\mathsf{V}}
\newcommand{\T}{\mathsf{T}}
\newcommand{\W}{\mathsf{W}}

\DeclareMathOperator{\A}{{\mathcal{A}}}

\DeclareMathOperator{\FP}{{FP}}

\newcommand{\inj}{\hookrightarrow}
\newcommand{\abs}[1]{\left| #1 \right|}

\newenvironment{romenum}
{\renewcommand{\theenumi}{\roman{enumi}}
\renewcommand{\labelenumi}{(\theenumi)}
\begin{enumerate}}{\end{enumerate}}

\newcounter{kcount}
\newenvironment{kenum}
{\begin{list}
{\textup{(K\thekcount)}}{\usecounter{kcount} 
\setlength{\leftmargin}{20pt}}}
{\end{list}}

\newcounter{qkcount}
\newenvironment{qkenum}
{\begin{list}
{\textup{(Q\theqkcount)}}{\usecounter{qkcount} 
\setlength{\leftmargin}{20pt}}}
{\end{list}}

\begin{document}

\title[Quasi-K\"{a}hler Bestvina-Brady groups]{%
Quasi-K\"{a}hler Bestvina-Brady groups}

\author[A. Dimca]{Alexandru Dimca}
\address{Laboratoire J.A.~Dieudonn\'{e}, UMR du CNRS 6621, 
Universit\'{e} de Nice--Sophia Antipolis, Parc Valrose,
06108 Nice Cedex 02, France}
\email{dimca@math.unice.fr}

\author[Stefan Papadima]{Stefan Papadima$^1$}
\address{Inst.~of Math.~Simion Stoilow, 
P.O. Box 1-764,
RO-014700 Bucharest, Romania}
\email{Stefan.Papadima@imar.ro}

\author[Alexander~I.~Suciu]{Alexander~I.~Suciu$^2$}
\address{Department of Mathematics,
Northeastern University,
Boston, MA 02115, USA}
\email{a.suciu@neu.edu}
\urladdr{http://www.math.neu.edu/\~{}suciu}

\thanks{$^1$Partially supported by CERES grant 4-147/12.11.2004
of the Romanian Ministry of Education and Research}

\thanks{$^2$Partially supported by NSF grant DMS-0311142}



\begin{abstract}
A finite simple graph $\G$ determines a right-angled Artin 
group $G_\G$, with one generator for each vertex $v$, 
and with one commutator relation $vw=wv$ for each 
pair of vertices joined by an edge.  The Bestvina-Brady 
group $N_\G$ is the kernel of the projection $G_\G\to \Z$, 
which sends each generator $v$ to $1$.  We establish  
precisely which graphs $\G$ give rise to quasi-K\"ahler 
(respectively, K\"ahler) groups $N_\G$. This yields 
examples of quasi-projective groups which are not 
commensurable (up to finite kernels) to the fundamental group of any 
aspherical, quasi-projective variety. 
\end{abstract}
\maketitle

\section{Introduction}
\label{sect:intro}

\subsection{}
\label{ss11}

Every finitely presented group $G$ is the fundamental group 
of a smooth, compact, connected manifold $M$ of dimension 
$4$ or higher.  Requiring that $M$ be a complex projective 
manifold (or, more generally, a compact K\"ahler manifold), 
puts extremely strong restrictions on what $G=\pi_1(M)$ 
can be; see \cite{ABCKT}.  Groups arising in this fashion 
are called projective (respectively, K\"ahler) groups. It is 
an open question whether K\"ahler groups need be projective. 

A related question, due to J.-P.~Serre, asks:  which finitely 
presented groups can be realized as fundamental groups of 
complements of normal crossing divisors in smooth, connected, 
complex projective varieties?  Groups arising in this fashion are 
called quasi-projective groups, while groups that can be realized 
as fundamental groups of complements of normal crossing 
divisors in compact K\"ahler manifolds are called quasi-K\"ahler 
groups. Again, it is an open question whether quasi-K\"ahler 
groups need be quasi-projective. 

Finally, we have the following question raised by J.~Koll\'ar 
in \cite{Ko95}, section 0.3.1:  given a (quasi-) projective 
group $G$, is there a group $\pi$, commensurable to $G$ 
up to finite kernels, and admitting a $K(\pi,1)$ which is a 
quasi-projective variety? (Two groups are 
commensurable up to finite kernels if there exists a zig-zag 
of homomorphisms, each with finite kernel and 
cofinite image, connecting one group to the other. 
See Remark \ref{rem:comm} for more details.)

\subsection{}
To a finite simple graph $\G=(\V,\E)$, there is associated a 
right-angled Artin group, $G_\G$, with a generator 
$v$ for each vertex $v\in \V$, and with a commutator 
relation $vw=wv$ for each edge $\{v,w\}\in \E$.  
The groups $G_\G$ interpolate between the free groups 
$F_n$ (corresponding to the discrete graphs on $n$ vertices, 
$\overline{K}_n$), and the free abelian groups $\Z^n$ 
(corresponding to the complete graphs on $n$ vertices, 
$K_n$).  These groups  are completely determined by their 
underlying graphs: $G_\G\cong G_{\G'}$ if and only if 
$\G\cong \G'$; see \cite{KMNR, Dr}. 

In a previous paper \cite{DPS}, we settled Serre's question 
for the class of right-angled Artin groups:  $G_\G$ is 
quasi-projective iff  $G_\G$ is quasi-K\"ahler iff  
$\G$ is a complete multipartite graph (i.e., a join 
$K_{n_1,\dots, n_r}=\overline{K}_{n_1} * \cdots * 
\overline{K}_{n_r}$ of discrete graphs). 
The object of this note is to answer Serre's question 
for the closely related class of Bestvina-Brady 
groups.  We also answer Koll\'ar's question, in 
the context of quasi-projective groups.

\subsection{}
Let  $\nu\colon G_{\G}\to \Z$ be the  homomorphism 
which sends each generator $v\in \V$ to $1$. 
The Bestvina-Brady group (or, Artin kernel) associated to 
$\G$ is defined as $N_{\G}=\ker (\nu)$. Unlike 
right-angled Artin groups, the groups $N_\G$ 
are {\em not}\/ classified by the underlying graphs $\G$. 
Moreover, these groups are complicated enough 
that a counterexample to either the Eilenberg-Ganea 
conjecture or the Whitehead conjecture can be 
constructed  from them, according to the fundamental 
paper of Bestvina and Brady \cite{BB}. 

Our answer to Serre's problem for this class of groups is 
as follows.

\begin{theorem}
\label{thm:bbqk} 
Let  $\Gamma$ be a finite simplicial graph, with associated 
Bestvina-Brady group $N_{\Gamma}$. The following are equivalent:
\begin{qkenum}
\item  \label{qk1} 
The group $N_{\Gamma}$ is quasi-K\"ahler.
\item  \label{qk2} 
The group $N_{\Gamma}$ is quasi-projective.
\item  \label{qk3} 
The graph  $\Gamma$ is either a tree, or a complete 
multipartite graph $K_{n_1,\dots ,n_r}$, with either some 
$n_i=1$, or all $n_i\ge 2$ and $r\ge 3$.
\end{qkenum}
\end{theorem}

Implication (Q\ref{qk2}) $\Rightarrow$ (Q\ref{qk1}) is 
clear. Implication (Q\ref{qk1}) $\Rightarrow$ (Q\ref{qk3})  
is established in Section~\ref{sect:proof}; the proof is based 
on certain cohomological obstructions to realizability 
by quasi-K\"ahler manifolds, developed in \cite{DPS}, 
and on computations of algebraic invariants for 
Bestvina-Brady groups, done in \cite{PSbb}.  
Implication (Q\ref{qk3}) $\Rightarrow$ (Q\ref{qk2}) 
is established in Section~\ref{sect:qp}; the proof is 
based on a result from \cite{D} on the topology of the 
mapping $f\colon\C^n \to \C$ induced by the defining 
equation of an affine hyperplane arrangement in $\C^n$.

\subsection{}
As a consequence, we can classify the 
Bestvina-Brady groups arising as fundamental groups 
of quasi-K\"ahler manifolds.  

\begin{corollary}
\label{cor:bbclass} 
The class of quasi-K\"ahler Bestvina-Brady groups equals the 
union of the following disjoint classes (which contain no repetitions):

\begin{enumerate}
\item  \label{b1} 
$\Z^{r}$, with $r\ge 0$;
\item  \label{b2} 
$\prod_{i=1}^{s} F_{n_i}$, with all $n_i>1$;
\item  \label{b3} 
$\Z^r\times \prod_{i=1}^{s} F_{n_i}$, with $r>0$ and all $n_i>1$;
\item  \label{b4} 
$N_{K_{n_1,\dots ,n_r}}$, with all $n_i\ge 2$ and $r\ge 3$.
\end{enumerate}
\end{corollary}

In the K\"ahler case, the above theorem takes the 
following, simpler form. 

\begin{corollary}
\label{cor:bbkahler} 
Let  $\Gamma$ be a finite simplicial graph, with associated 
Bestvina-Brady group $N_{\Gamma}$. The following are equivalent:
\begin{kenum}
\item  \label{k1} 
The group $N_{\Gamma}$ is K\"ahler.
\item  \label{k1.5} 
The group $N_{\Gamma}$ is projective.
\item  \label{k2} 
The graph  $\Gamma$ is a complete 
graph $K_{n}$ with $n$ odd.
\item  \label{k3} 
The group $N_{\Gamma}$ is free abelian of even rank.
\end{kenum}
\end{corollary}

Our answer to Koll\'ar's question, in the context of 
quasi-projective Bestvina-Brady groups, is as follows.  

\begin{corollary}
\label{cor:bbkpi1} 
Let $N_{\Gamma}$ be a quasi-projective (equivalently, 
quasi-K\"ahler)  Bestvina-Brady group. Then, in terms 
of the classification from Corollary \ref{cor:bbclass}:  
\begin{romenum}
\item \label{i}
$N_{\Gamma}$ is of type \eqref{b1}, \eqref{b2}, 
or \eqref{b3} if and only if it is the fundamental group 
of an aspherical, smooth quasi-projective variety.

\item  \label{ii}
$N_{\Gamma}$ is of type \eqref{b4} if and only if 
it is not commensurable up to finite kernels to the 
fundamental group of any aspherical, quasi-projective variety.
\end{romenum}
\end{corollary}

The simplest example of a group of type \eqref{ii} is 
the group $G=N_{K_{2,2,2}}$, already studied by Stallings in 
\cite{St}.  As noted in \cite{PSdecomp}, this group is the 
fundamental group of the complement of an arrangement 
of  lines in $\mathbb{CP}^2$; thus, $G$ is quasi-projective.  
On the other hand, as shown by Stallings, $H_3(G,\Z)$ is 
not finitely generated; thus, there is no $K(G,1)$ space 
with finite $3$-skeleton.

The three corollaries above are proved at the end of 
Section \ref{sect:proof}.

\section{Resonance obstructions}
\label{sect:review}

We start by reviewing the cohomological obstructions to 
realizability by quasi-K\"ahler manifolds that we will need 
in the sequel.  

Let $G$ be a finitely presented group. 
 Denote by $A=H^*(G,\C)$ the cohomology algebra of $G$. 
For each $a\in A^1$, we have $a^2=0$, and so 
right-multiplication by $a$ defines a cochain complex
$(A,a)\colon A^0 \xrightarrow{a}  A^1\xrightarrow{a} A^2$. 
Let $\RR_1(G)$ be the set of points $a\in A^1$ where 
this complex fails to be exact, 
\begin{equation*} 
\label{eq:rv}
\RR_1(G)=\{a \in A^1 \mid H^1(A, a) \ne 0\}.
\end{equation*}
The set $\RR_1(G)$ is a homogeneous algebraic 
variety in the affine space $A^1=H^1(G,\C)$, called 
the {\em (first) resonance variety} of $G$.

The group $G$ is said to be {\em $1$-formal} if its Malcev 
Lie algebra is quadratically presented.  This important notion 
goes back to the work of Quillen \cite{Q} and Sullivan \cite{S} 
on rational homotopy theory; for a comprehensive account 
and recent results on $1$-formal groups, see \cite{PS04}.

The best-known examples of $1$-formal groups are the 
K\"ahler groups, see \cite{DGMS}. Other classes of examples 
include Artin groups \cite{KM} and finitely-presented 
Bestvina-Brady groups \cite{PSbb}.  The $1$-formality 
property is preserved under free products and direct products, 
see \cite{DPS}. 

\begin{theorem}[\cite{DPS}]
\label{thm:posobstr} 
Let $M$ be a quasi-K\"ahler manifold. 
Set $G=\pi_1(M)$ and let 
$\RR_1(G)=\bigcup\nolimits_{\alpha}\RR^{\alpha}$
be the decomposition of the first resonance variety 
of $G$ into irreducible components.  If $G$ is $1$-formal, 
then the following hold.

\begin{enumerate}
\item  \label{a1}  
Every positive-dimensional component $\RR^{\alpha}$ 
is a $p$-isotropic linear subspace of $H^1(G, \C)$, 
of dimension at least $2p+2$, for some 
$p=p(\alpha) \in \{0,1\}$. 

\item  \label{a2} 
If $M$ is K\"ahler, then only $1$-isotropic components 
can occur.
\end{enumerate}
\end{theorem}

Here, we say that a nonzero subspace $V\subseteq  H^1(G,\C)$ 
is $0$- (respectively, $1$-) {\em isotropic} if the 
restriction of the cup-product map, 
$\cup_G\colon \bigwedge^2 V\to \cup_G (\bigwedge^2 V)$, 
is equivalent to 
$\cup_C \colon \bigwedge^2 H^1(C,\C)\to H^2(C,\C)$,
where $C$ is a non-compact (respectively, compact) smooth, 
connected complex curve.

\section{Cohomology ring and resonance 
varieties for $G_\G$ and $N_\G$}
\label{sect:review2}

Let $\G=(\V,\E)$ be a finite simple graph, 
with $G_\G=\langle v\in \V \mid vw=wv 
\text { if $\{v,w\}\in \E$}\rangle$ 
the corresponding right-angled Artin group, 
and $N_\G=\ker (\nu \colon G_\G \to \Z)$ the 
corresponding Bestvina-Brady group.  In this 
section, we review some facts about the 
cohomology in low degrees and the resonance 
varieties for the groups $G_\G$ and $N_\G$.

\subsection{} 
\label{subsec:resartin}
 
Write $H_{\V}=H^1(G_\G, \C)$.  
If $\W$ is a subset of $\V$, write $\G_\W$ for the induced 
subgraph on vertex set $\W$, and $H_{\W}$ 
for the corresponding coordinate subspace of $H_{\V}$. 

\begin{theorem}[\cite{PS06}]
\label{thm:res artin}
Let $\G$ be a finite graph.  Then:
\begin{equation*} 
\label{eq:res artin}
\RR_1(G_{\G}) = \bigcup_{\stackrel{\W\subset \V}{
\G_{\W} \textup{ disconnected}}} H_{\W}.
\end{equation*} 
\end{theorem}

It follows that the irreducible components of $\RR_1(G_\G)$ 
are indexed by the subsets $\W\subset \V$, 
maximal among those for which $\G_{\W}$ is disconnected. 

Using Theorems \ref{thm:posobstr} and \ref{thm:res artin}, 
we gave in \cite{DPS} the following answer to Serre's problem 
for right-angled Artin groups. 

\begin{theorem}[\cite{DPS}] 
\label{thm:artinserre} 
Let $\Gamma=(\V,\E)$ be a finite simplicial graph, with 
associated right-angled Artin group $G_{\Gamma}$. The 
following are equivalent.

\begin{romenum}
\item  \label{s1} 
The group $G_{\Gamma}$ is quasi-K\"ahler.
\item  \label{s2} 
The group $G_{\Gamma}$ is quasi-projective.
\item  \label{sob}  
The group $G_{\Gamma}$ satisfies the isotropicity condition 
from Theorem \ref{thm:posobstr}\eqref{a1}. 
\item  \label{s3}  
The group $G_{\Gamma}$ is a product of finitely generated 
free groups.
\item  \label{s4} 
The graph  $\Gamma$ is a complete multipartite graph 
$K_{n_1,\dots ,n_r}=\overline{K}_{n_1}*\dots *\overline{K}_{n_r}$.
\end{romenum}
\end{theorem}

\subsection{} 
\label{subsec:resbb}

In order to apply Theorem \ref{thm:posobstr} to the 
Bestvina-Brady groups, we need a description 
of their resonance varieties. 

The inclusion $\iota\colon N_{\G} \to G_{\G}$ induces 
a homomorphism 
$\iota^*\colon H^1(G_{\G},\C) \to  H^1(N_{\G},\C)$.
Note that if $\abs{\V}=1$, then $N_{\G}=\{1\}$ and thus 
$\RR_1(N_{\G})$ is empty.

\begin{theorem}[\cite{PSbb}] 
\label{thm:resbb}
Let $\G$ be a graph, with connectivity $\kappa(\G)$. 
Suppose the flag complex $\Delta_{\G}$ is simply-connected, 
and $\abs{\V}>1$. 
\begin{enumerate}
\item \label{r1}
If $\kappa(\G)=1$, then $\RR_1(N_{\G})=H^1(N_{\G},\C)$.  
\item \label{r2}
If $\kappa(\G)>1$, then the irreducible components of 
$\RR_1(N_{\G})$ are the subspaces $H'_{\W}=\iota^*(H_{\W})$, 
of dimension $\abs{\W}$, one for each subset 
$\W\subset \V$, maximal among those for 
which $\Gamma_{\W}$ is disconnected.  
\end{enumerate}
\end{theorem}

\subsection{} 
\label{subsec:cohoartinbb}

To analyze the isotropicity properties of the 
resonance components in the above theorem, 
we need a precise description 
of the cohomology ring in degrees $1$ and $2$, for both 
right-angled Artin groups and Bestvina-Brady groups. 

For $G_\G$, the cohomology 
ring can be identified with the exterior Stanley-Reisner 
ring of the flag complex:  $H^*(G_\G)$ is the quotient 
of the exterior algebra on generators $v^*$ in degree $1$, 
indexed by the vertices $v\in \V$, modulo the ideal generated 
by the monomials $v^* w^*$ for which $\{v,w\}$ is not an 
edge of $\G$.  It follows that an additive basis for $H^k(G_\G)$ 
is indexed by the complete $k$-subgraphs of $\G$. 

Viewing the homomorphism $\nu\colon G_\G \to \Z\inj \C$ 
as an element in $H^1(G_\G,\C)$, we have the following. 

\begin{theorem}[\cite{PSbb}] 
\label{thm:cohobb}
Suppose $\pi_1(\Delta_{\G})=0$. Then 
$\iota^*\colon H^*(G_\G,\C)\to H^*(N_\G,\C)$ induces 
a ring homomorphism $\iota^*\colon 
H^*(G_\G,\C)/(\nu\cdot H^*(G_\G,\C))\to H^*(N_\G,\C)$, 
which is an isomorphism in degrees $*\le 2$.
\end{theorem}
In other words, $\iota^*\colon H^1(G_\G,\C)\to H^1(N_\G,\C)$ 
is onto, has one dimensional kernel generated by $\nu$, 
and fits into the following commuting diagram:
\begin{equation}
 \label{eq:cupgn}
 \xymatrix{
\bigwedge^2 H^1(G_\G,\C) 
 \ar@{->>}[r]^{\wedge^2 \iota^*} 
\ar[d]^{\cup_{G_\G}} & 
 \bigwedge^2 H^1(N_\G,\C) 
 \ar[d]^{\cup_{N_\G}}\\
   H^2(G_\G,\C) \ar@{->>}[r]
   &  H^2(G_\G,\C) /\nu H^1(G_{\G},\C)\cong  H^2(N_\G,\C)
 }
\end{equation}

\section{The quasi-K\"ahler condition for the groups $N_\G$}
\label{sect:proof}

In this section, we establish implication 
(Q\ref{qk1}) $\Rightarrow$ (Q\ref{qk3}) 
from Theorem \ref{thm:bbqk} in the Introduction: 
the implication follows at once from Lemmas \ref{lem:dic1} 
and \ref{lem:dic3} below.  We conclude the section with 
proofs of Corollaries \ref{cor:bbclass}, \ref{cor:bbkahler}, 
and \ref{cor:bbkpi1}. 

\subsection{}
Let $\Gamma=(\V,\E)$ be a finite simple graph.  
Denote by $\T$ the set of triangles in $\G$. 
Suppose the Bestvina-Brady group $N_\G$ is 
finitely presented.  The following then hold:

\subsubsection*{\textup{I}}
The flag complex $\Delta_\G$ 
must be simply-connected, as shown in \cite{BB}. 

\subsubsection*{\textup{II}}
An explicit finite presentation for $N_\G$ is given 
by Dicks and Leary \cite{DL}. Fix a linear order on 
$\V$, and orient the edges increasingly.  
A triple of edges $(e,f,g)$ forms a directed triangle 
if $e=\{u,v\}$, $f=\{v,w\}$, $g=\{u,w\}$, and  $u<v<w$.
Then:
\begin{equation}
\label{eq:bb group}
N_{\G}= \langle e \in \E \mid 
ef=fe, \: ef=g\: \text{ if $(e,f,g)$ is a directed triangle}\,
\rangle.
\end{equation}

\subsubsection*{\textup{III}}
The group $N_\G$ is $1$-formal, as shown in 
Proposition 6.1 from \cite{PSbb}. 

\subsection{}
Now suppose $N_\G=\pi_1(M)$, with $M$ a 
quasi-K\"ahler manifold.  Then $N_\G$ is 
finitely presented, and thus $1$-formal, 
by the above. Therefore, Theorem \ref{thm:posobstr} 
applies to $N_\G$.

\begin{lemma}
\label{lem:dic1}
Suppose $N_\G$ is a quasi-K\"ahler group. 
If there exists an edge $e_0\in \E$ not included in 
any triangle in $\T$, then $\G$ is a tree.
\end{lemma}

\begin{proof}
We may safely assume $\G$ is not $K_1$ or $K_2$, 
to avoid degenerate cases.
Since $\pi_1(\Delta_\G)=0$, it is enough to show that 
$\G$ has no triangles. Indeed, if $\T=\emptyset$, then 
the $2$-skeleton of $\Delta_\G$ coincides with $\G$. 
Hence, $\pi_1(\G)=0$, and so $\G$ is a tree.

Write $\E'=\E\setminus \{e_0\}$. 
From the Dicks-Leary presentation \eqref{eq:bb group}, 
we find that $N_{\G}=\Z * N' $, where $\Z$ is the cyclic 
group generated by $e_0$, and 
\begin{equation*}
\label{eq:nprime}
N'=\langle e\in \E' \mid 
ef=fe, \: ef=g\: \text{ if $(e,f,g)$ is a directed triangle}\,
\rangle.
\end{equation*}

We may assume that $b_1(N')\ne 0$, for otherwise 
$\Gamma =K_2$. From \cite[Corollary~5.4]{PS06}, 
we find that  $\RR_1(N_\G)=H^1(N_\G, \C)$. Moreover, 
by \cite[Lemma~7.4]{DPS}, the cup-product map 
$\cup_{N_\G}$ vanishes. 

Suppose $\tau=\{u,v,w\}$ is a triangle in $\T$. 
By the above,  $\iota^*(u^*) \cup_{N_\G} \iota^*(v^*) = 0$ 
in $H^2(N_\G, \C)$. 
Write $\nu=\sum_{z\in \V} z^*\in H^1(G_\G,\C)$.
From diagram \eqref{eq:cupgn}, we find that 
$u^*\cup_{G_\G} v^* \cup_{G_\G} \nu =0$ in 
$H^3(G_\G, \C)$.  Hence:
\[
\sum_{z\in \V} u^*\cup_{G_\G} v^* \cup_{G_\G} z^* = 0 . 
\]
But this sum contains the basis element $u^*v^*w^*$ 
of $H^3(G_\G, \C)$, a contradiction.
\end{proof}

\begin{lemma}
\label{lem:dic2}
Suppose $N_\G$ is a quasi-K\"ahler group, 
and every edge of $\G$ belongs to a triangle. 
Then $\kappa(\G)>1$.
\end{lemma}

\begin{proof}
Suppose $\G$ has connectivity $1$.  Then, by 
Theorem \ref{thm:resbb}\eqref{r1}, we must have 
$\RR_1(N_{\G})=H^1(N_{\G},\C)$. Furthermore, this linear 
space must be either $0$- or $1$-isotropic, by Theorem 
\ref{thm:posobstr}\eqref{a1}.  That is, the cup-product map 
$\cup_{N_\G}\colon \bigwedge^2 H^1(N_\G,\C) \to H^2(N_\G,\C)$ 
either vanishes, or has $1$-dimensional image and is 
non-degenerate. 

On the other hand, we know from Theorem \ref{thm:cohobb} that 
the map $\cup_{N_\G}$ is surjective. Thus, we must have 
$b_2(N_\G)\le 1$. 
Moreover, by Proposition 7.1 in \cite{PSbb}, we have 
$b_2(N_\G)= \abs{\E}-\abs{\V}+1$, since $\Delta_\G$ 
is simply-connected. Thus, $ \abs{\V} -1\le \abs{\E}\le \abs{\V}$. 

If $\abs{\E}=\abs{\V}- 1$, then $\G$ is a tree, contradicting 
the assumption that every edge belongs to a  triangle.  
If $\abs{\E}=\abs{\V}$, then $\G$ is obtained 
from a tree by adding exactly one edge. 
Using again the assumption that every edge 
belongs to a  triangle, it is readily seen that 
$\G=K_3$.  But $\kappa(K_3)>1$, a contradiction. 
\end{proof}

\begin{lemma}
\label{lem:dic3}
Suppose $N_\G$ is a quasi-K\"ahler group, 
and every edge of $\G$ belongs to a triangle. Then
$\G=K_{n_1,\dots ,n_r}$, with either some 
$n_i=1$, or all $n_i\ge 2$ and $r\ge 3$.
\end{lemma}

\begin{proof}
If $\G$ is a complete graph $K_n=K_{1,\dots,1}$, 
there is nothing to prove. Otherwise, all irreducible 
components $H'_\W$ of $\RR_1(N_\G)$ are positive-dimensional, 
by Lemma \ref{lem:dic2} and Theorem \ref{thm:resbb}\eqref{r2}.  
Moreover, each subspace $H'_\W$ is $0$- or $1$-isotropic, 
by Theorem \ref{thm:posobstr}\eqref{a1}. 

\begin{claim}
The resonance obstruction from Theorem 
\ref{thm:posobstr}\eqref{a1} 
holds for $G_\G$. 
\end{claim}

Assuming this claim, the proof is completed as follows. 
By Theorem \ref{thm:artinserre}, the graph $\G$ 
must be a complete multipartite graph 
$K_{n_1,\dots, n_r}$. If all $n_i\ge 2$, then 
necessarily $r\ge 3$, by the simply-connectivity 
assumption on $\Delta_\G$, and we are done. 

Let $\bigcup_{\W} H_\W$ be the decomposition into irreducible 
components of $\RR_1(G_\G)$ from Theorem \ref{thm:res artin}, 
where $\W$ runs through the non-empty, maximal subsets of 
$\V$ for which $\G_\W$ is disconnected. We will prove that 
each induced subgraph $\Gamma_\W$ is a discrete graph.  
This will imply that each subspace $H_\W$ is a $0$-isotropic 
subspace of $H^1(G_\G,\C)$, of dimension at least $2$,  
thereby showing that $G_\G$ satisfies the resonance 
obstructions, as claimed. 

Let $\G_\W=\bigsqcup_{j} \G_{\W_j}$ be the decomposition into 
connected components of the disconnected graph $\G_\W$. 
Denote by $m(\W)$ the number of non-discrete connected components. 

First suppose $m(\W)\ge 2$. Then, there exist disjoint edges, 
$e_1$ and $e_2$, in $\G_\W$. Recall that each edge $e\in \E$ 
can be viewed as a basis element of $H^2(G_\G,\C)$. From the 
isotropicity of $H'_{\W}$, we deduce that $\iota^*(e_1)$ and 
$\iota^*(e_2)$ are linearly dependent elements in $H^2(N_\G,\C)$.  
But this is impossible.  Indeed, if 
$\alpha_1 \iota^*(e_1) + \alpha_2 \iota^*(e_2) = 0$ is a 
linear dependence, then it follows from diagram \eqref{eq:cupgn} 
that $(\alpha_1 e_1 + \alpha_2 e_2)\cdot \nu = 0$ 
in $H^3(G_\G,\C)$. Expanding the product, we find that
\begin{equation}
\label{eq:lindep}
\alpha_1 \sum_{\tau \in \T \colon \tau \supset e_1} \tau + 
\alpha_2 \sum_{\tau \in \T \colon \tau \supset e_2} \tau = 0\, .
\end{equation}
On one hand, we know that there are 
triangles $\tau_1$ and $\tau_2$ containing edges $e_1$ and $e_2$, 
respectively.  On the other hand, all basis elements $\tau$ 
appearing in \eqref{eq:lindep} are distinct, since 
$e_1\cap e_2=\emptyset$. Thus,  $\alpha_1=\alpha_2=0$. 

Now suppose $m(\W)=1$.  Write $\W=\W'\bigsqcup \W''$, 
with $\G_{\W'}$ discrete and non-empty, 
$\G_{\W''}$ containing at least an edge, 
and with no edge joining $\W'$ to $\W''$. 
We then  have a non-trivial orthogonal decomposition 
with respect to $\cup_{N_\G}$, 
\begin{equation}
\label{eq:ortho}
H'_{\W}=H'_{\W'}\oplus H'_{\W''}. 
\end{equation}
As above, we find that the restriction of $\cup_{N_\G}$ to 
$\bigwedge^2 H'_{\W''}$ is non-zero.  Thus, $H'_{\W}$ 
cannot be $0$-isotropic. On the other hand, any non-zero 
element in $H'_{\W'}$ is orthogonal to all of $H'_{\W}$; 
thus the restriction of $\cup_{N_\G}$ to $H'_{\W}$ is 
a degenerate form. Hence, $H'_{\W}$ cannot be $1$-isotropic, 
either. This contradicts the isotropicity property of $N_\G$. 

Thus, $m(\W)=0$, and $\G_\W$ is discrete, as asserted. 
This ends the proof. 
\end{proof}

\subsection{Proof of Corollary \ref{cor:bbclass}}
\label{subs:corbbclass}

The classification of quasi-K\"ahler Bestvina-Brady groups 
follows from Theorem \ref{thm:bbqk}.  The passage from
graphs to groups is achieved by using the following simple 
facts: if $\G$ is a tree on $n$ vertices, then $N_{\G}=F_{n-1}$;
if $\G =\G_1 *\G_2$, 
then $G_{\G}=G_{\G_1} \times G_{\G_2}$;
if $\G= K_1 * \G'$, then $N_{\G}=G_{\G'}$.

Next, note that the groups $G$ in classes \eqref{b1}--\eqref{b3} 
have a finite $K(G,1)$, while those in class  \eqref{b4} don't, 
by \cite{BB}. The groups in \eqref{b1}--\eqref{b3} 
are distinguished by their Poincar\'e polynomials. 
As for the groups in class  \eqref{b4}, they are 
distinguished by their resonance varieties. Indeed, if 
$\G=K_{n_1,\dots,n_r}$, with all $n_i\ge 2$ and $r\ge 3$, 
then  $\kappa(\G)>1$ and $\pi_1(\Delta_\G)=0$; 
thus, by Theorem \ref{thm:resbb}\eqref{r2}, the  
variety  $\RR_1(N_\G)$ decomposes into $r$ 
irreducible components, of dimensions $n_1,\dots,n_r$. 

This finishes the proof of Corollary \ref{cor:bbclass}. 
\hfill\qed

\subsection{Proof of Corollary \ref{cor:bbkahler}}
\label{subs:corbbk}

Implications (K\ref{k2}) $\Rightarrow$ (K\ref{k3})
and  (K\ref{k1.5}) $\Rightarrow$ (K\ref{k1}) are clear. 
The fact that the torus $T^{2n}$ is a smooth projective 
manifold proves (K\ref{k3}) $\Rightarrow$ (K\ref{k1.5}).  

We are left with proving (K\ref{k1}) $\Rightarrow$ (K\ref{k2}). 
Suppose $N_\G$ is a K\"ahler group. Then of course $N_\G$ 
is quasi-K\"ahler, so $\G$ must be one of the graphs described
in Theorem \ref{thm:bbqk} (Q\ref{qk3}).

Assume $\G$ is neither a complete graph,
nor a complete multipartite graph $K_{n_1,\dots,n_r}$ 
with all $n_i\ge 2$ and $r\ge 3$. By the same argument 
as in the proof of Corollary \ref{cor:bbclass}, we infer 
that $N_{\G}$ must be of the form $F_n\times N'$, with 
$n\ge 2$. But this cannot be a K\"ahler group, by a 
result of Johnson and Rees, see \cite[Theorem 3]{JR}.

Suppose now $\G=K_{n_1,\dots,n_r}$, with all $n_i\ge 2$ 
and $r\ge 3$.  It follows from Theorem \ref{thm:posobstr}\eqref{a2} 
that all positive-dimensional irreducible components of 
$\RR_1(N_\G)$ must be $1$-isotropic. But this contradicts 
Theorem \ref{thm:resbb}\eqref{r2}, which predicts only 
$0$-isotropic components.

Finally, if $\G= K_n$, note that $b_1(N_\G)=n-1$.
Since the odd Betti numbers of a compact K\"ahler 
manifold are even, $n$ must be odd. 

This finishes the proof of Corollary \ref{cor:bbkahler}. 
\hfill\qed

\subsection{Proof of Corollary \ref{cor:bbkpi1}}  
Let $G=N_\Gamma$ be a quasi-projective Bestvina-Brady 
group, as classified in Corollary \ref{cor:bbclass}.  

If $G$ is of type \eqref{b1}, \eqref{b2}, or \eqref{b3}, then 
$G$ is the fundamental group of a space of the form 
$M=\prod_{i=1}^{q} (\C\setminus \{\text{$n_i$ points}\})$, 
an aspherical, smooth quasi-projective variety.  

If $G$ is of type \eqref{b4}, then $G$ is not of type 
$\FP_{\infty}$, by  \cite{BB} (see also Proposition~7.1 
and Remark~5.8 in \cite{PSbb} for a direct proof).  
As is well-known, a finite-index subgroup $\pi'$ of a 
group $\pi$ is of type $\FP_{\infty}$ if and only if  $\pi$ 
is of type $\FP_{\infty}$; see \cite[Prop.~VIII.5.1]{Br}. 
In particular, finite groups are of type $\FP_{\infty}$.
Moreover, if $\pi$ is an extension of $\pi'$ by a finite group,
then  $\pi$ is of type $\FP_{\infty}$ if and only if  $\pi'$ is;
see \cite[Proposition 2.7]{Bi}.

Thus, $G$ cannot be commensurable up to finite kernels to a group 
$\pi$ of type $\FP_{\infty}$, and so, a fortiori, to a 
group $\pi$ admitting a finite-type $K(\pi,1)$. 
We conclude that $G$ is not commensurable up to finite kernels to 
any group of the form $\pi=\pi_1(M)$, where 
$M$ is an aspherical, quasi-projective variety.

This finishes the proof of Corollary \ref{cor:bbkpi1}.  
\hfill\qed

\begin{remark}
\label{rem:comm}
Two groups, $G$ and $\pi$, are 
{\em elementarily commensurable up to finite kernels}\/ 
if there is a homomorphism $\varphi\colon G\to \pi$ with 
finite kernel and with finite-index image; the associated 
equivalence relation is the one defined at the end of 
\S\ref{ss11}.  A stronger equivalence relation is defined 
as follows: $G$ and $\pi$ are {\em commensurable}\/ 
if there exist finite-index subgroups, $G'<G$ and $\pi'<\pi$, 
such that $G'\cong \pi'$. 

Being elementarily commensurable up to finite kernels 
is not a symmetric relation.  An example, due to Milnor, 
can be found within the class of finitely generated, 
nilpotent, torsion-free groups; see \cite[p.~25]{HMR}. 
It is worth noting that, in this class of groups, the two  
commensurability relations coincide: they both amount to 
saying that $G$ and $\pi$ have the same rationalization; 
see \cite[Theorem I.3.3]{HMR} and \cite[Lemma 2.8 on p.~19]{Ba}. 
\end{remark}

\section{Realizability by quasi-projective varieties}
\label{sect:qp}

In this section, we prove implication (Q\ref{qk3}) $\Rightarrow$ 
(Q\ref{qk2}) from Theorem \ref{thm:bbqk}.  The implication 
is easily proved for groups belonging to the classes 
\eqref{b1}--\eqref{b3} listed in Corollary \ref{cor:bbclass}, 
by taking products of punctured projective lines. To realize 
the groups in class \eqref{b4} by smooth quasi-projective 
varieties, we will appeal to results on the topology of
fibers of polynomial mappings, see \cite{B, D}. 

Let $\A=\{H_1,\dots ,H_d\}$ be an affine essential hyperplane 
arrangement in $\C^{n}$.  We set $M=\C^{n} \setminus X$, 
with $X$ being the union of all the hyperplanes in $\A$.
Consider a $d$-tuple of positive integers, $e=(e_1,\dots ,e_d)$, 
with $\gcd (e_1,\dots ,e_d)=1$. 
For any $i=1,\dots ,d$, let $\ell _i=0$ be an equation 
for the hyperplane $H_i$ and consider the product 
\[
f_e=\prod _{i=1}^{d}\ell_i^{e_i} \in \C[x_1,\dots ,x_n].
\]

Let $d_e=e_1+\cdots +e_d$ be the degree of the polynomial $f_e$. 
Note that  $f_e=0$ is a (possibly non-reduced) equation for the 
union $X$.  The induced morphism 
\[
\nu_e=(f_e)_{\sharp}\colon \pi _1(M) \to \pi _1(\C^*)=\Z
\]
sends an elementary (oriented) loop about the hyperplane $H_j$ 
to $e_j \in \Z$.

When the arrangement $\A$ is central, i.e., $0 \in H_i$ 
for all $i=1,\dots ,d$, the above polynomial $f_e$ is 
homogeneous, and there is a lot of interest in the 
associated Milnor fiber $F_e=f_e^{-1}(1)$. 
In particular, we have the following exact sequence of groups:
\[
 \xymatrix{1 \ar[r]&  \pi _1(F_e) \ar[r]&  \pi _1(M)\ar[r]& 
  \pi _1(\C^*)=\Z \ar[r]&  0}.
 \]
The main result of this section is the following. 

\begin{theorem} 
\label{thm:mf}
For any affine essential hyperplane arrangement in
$\C^{n}$ with $n \geq 3$, any relatively prime integers 
$e$ as above and any $t \in \C^*$, there is an
exact sequence of groups
\[
 \xymatrix{1 \ar[r] &  \pi _1(F_{e,t}) \ar[r] &  \pi _1(M) 
 \ar[r]^(.43){\nu_e}
 &   \pi _1(\C^*)=\Z \ar[r]&  0},
 \]
where $F_{e,t}=f_e^{-1}(t)$.
\end{theorem}

\begin{proof}
 It is known that the mapping $f_e\colon M \to \C^*$ has 
 only finitely many singular points, see \cite[Corollary 2.2]{D}.
Let $R>0$ be a real number such 
that the disc
\[
D_R=\{z\in \C~~;~~|z|\leq R\}
\]
contains in its interior all the critical values of $f_e$. 
We set $S^1_R=\{z\in \C~~;~~|z|= R\}$.
Let $M_R=M\cap f_e^{-1}(D_R)$ and $E_R=f_e^{-1}(S^1_R)$. 
The main result in \cite{D}, Theorem 2.1(ii),
says roughly that our polynomial mapping $f_e$ behaves 
like a proper mapping. This comes from the fact that all the 
fibers of $f_e$ are transversal to large enough spheres 
centered at the origin of $\C^n$, similarly to the case of 
tame polynomials considered in \cite{B}. In particular, 
this shows that the bifurcation set of the mapping 
$f_e\colon \C^n \to \C$, i.e., the minimal set to be 
deleted in order to have a locally trivial fibration, 
coincides with the set of critical values of $f_e$.
Specifically, the following properties hold.

\smallskip\noindent
(A)  The restriction of $f_e$ induces a locally 
trivial fibration $f_e\colon E_R \to S^1_R$.
In particular, we get an exact sequence
\[
 \xymatrix{1 \ar[r]&  \pi _1(F_{e,R}) \ar[r]&  \pi _1(E_R) 
 \ar[r]& \pi _1(S^1_R)=\Z  \ar[r]&  0}.
\]

\smallskip\noindent
(B)\,  The inclusion $M_R \hookrightarrow M$ is 
a homotopy equivalence.

\smallskip\noindent
(C)  The space $M_R$ has the homotopy type of a 
CW-complex obtained from $E_R$ by attaching
finitely many $n$-cells. This is a general trick in this situation, 
see \cite{B}. In particular, the inclusion $E_R  \hookrightarrow  M$ 
induces an isomorphism at the level of fundamental groups if $n\geq 3$.

\smallskip\noindent
(D)  Any smooth fiber $F_{e,t}=f_e^{-1}(t)$ of $f_e$ is 
diffeomorphic to the fiber $F_{e,R}=f_e^{-1}(R)$. And any 
singular fiber $F_{e,t}=f_e^{-1}(t)$ of $f_e$ for $t \ne 0$ is 
obtained from the fiber $F_{e,R}=f_e^{-1}(R)$ by adding 
a finite number of $n$-cells. In particular, if $n\geq 3$, there are
natural isomorphisms $\pi _1(F_{e,t}) \to \pi _1(F_{e,R})$ for 
any $t \in \C^*$.

\smallskip
These properties allow us to transform the exact sequence 
in (A) to get our claimed exact sequence.
\end{proof}

The quasi-projectivity of the groups $N_{K_{n_1,\dots,n_r}}$  
($r\ge 3$) follows via a simple remark:  groups of the form 
$F_{n_1}\times \cdots \times F_{n_r}$ can be realized by affine 
arrangements; see the proof of the Corollary below.  
Related results, on the topology of {\em infinite}\/ 
cyclic covers for arrangements having fundamental 
groups of this form, may be found in \cite{CDP}.

\begin{corollary}
\label{cor=4qproj}
The Bestvina-Brady groups associated to the graphs 
$K_{n_1,\dots,n_r}$ with $r\ge 3$ are quasi-projective.
\end{corollary}

\begin{proof} 
Let $\A$ be the arrangement in $\C^r$ defined by the polynomial
\[
f = (x_1 -1)\cdots (x_1-n_1)\cdots (x_r -1)\cdots (x_r-n_r) .
\]
Theorem \ref{thm:mf} may be applied to $\A$ and $e=(1,\dots, 1)$.
Clearly, $\pi_1(M)$ may be identified with the right-angled 
Artin group $G_\G=F_{n_1}\times \cdots \times F_{n_r}$ 
defined by  the graph $\G= K_{n_1,\dots,n_r}$, 
in such a way that  $\nu_e$ becomes identified with 
the length homomorphism $\nu \colon G_\G \to \Z$. 
Hence, $N_\G$ is isomorphic to the fundamental group 
of the smooth fiber $F_{e,t}=f_e^{-1}(t)$ of $f_e$.
\end{proof}

This Corollary completes the proof of Theorem \ref{thm:bbqk}.

\begin{ack}
This work was done while the second author was 
visiting Northeastern University, in Spring, 2006. 
He thanks the Northeastern Mathematics Department 
for its support and hospitality during his stay. 

We thank J\'anos Koll\'ar for bringing up \cite{Ko95} to our 
attention, and for indicating its relevance to our work. 
\end{ack}

\end{document}